\documentclass[a4paper,11pt]{amsart}
\usepackage{graphicx}
\usepackage{mathptmx}
\usepackage{amsmath}
\usepackage{amssymb}
\usepackage{enumitem}
\usepackage{xcolor}

\newmuskip\pFqmuskip

\newcommand*\pFq[6][8]{%
  \begingroup 
  \pFqmuskip=#1mu\relax
  \mathcode`=\string"8000
  \begingroup\lccode`\~=`\,
  \lowercase{\endgroup\let~}\pFqcomma
  F^{#2}_{#3}{\left(\genfrac..{0pt}{}{#4}{#5}\bigg|#6\right)}%
  \endgroup
}
\newcommand{\pFqcomma}{\mskip\pFqmuskip}

\newtheorem{theorem}{Theorem}

\newtheorem{corollary}[theorem]{Corollary}

\newtheorem{remark}[theorem]{Remark}

\begin{document}

\title[]{New approach to $\lambda$-Stirling numbers}

\author{Dae San Kim}
\address{Department of Mathematics, Sogang University, Seoul 121-742, Republic of Korea}
\email{dskim@sogang.ac.kr}

\author{Hye Kyung Kim$^{2,*}$}
\address{Department of Mathematics Education, Daegu Catholic University, Gyeongsan 38430, Republic of Korea}
\email{hkkim@cu.ac.kr}

\author{Taekyun Kim $^{3,*}$}
\address{Department of Mathematics, Kwangwoon University, Seoul 139-701, Republic of Korea}
\email{tkkim@kw.ac.kr}

\thanks{* is corresponding author}

\begin{abstract}
The aim of this paper is to study the $\lambda$-Stirling numbers of both kinds which are $\lambda$-analogues of Stirling numbers of both kinds. Those numbers have nice combinatorial interpretations when $\lambda$ are positive integers. If $\lambda=1$, then the $\lambda$-Stirling numbers of both kinds reduce to the Stirling numbers of both kinds. We derive new types of generating functions of the $\lambda$-Stirling numbers of both kinds which are related to the reciprocals of the generalized rising factorials. Furthermore, some related identities are also derived from those generating functions. In addition, all the corresponding results to the $\lambda$-Stirling numbers of both kinds are obtained also for the $\lambda$-analogues of $r$-Stirling numbers of both kinds which are generalizations of those numbers.
\end{abstract}

\maketitle
 
\markboth{\centerline{\scriptsize New approach to $\lambda$-Stirling numbers}}
{\centerline{\scriptsize   D. S. Kim, H. K. Kim and T. Kim}}

\section{Introduction}
It is remarkable that explorations for degenerate versions of quite a few special numbers and polynomials, which began with the study of degenerate Bernoulli and degenerate Euler polynomials by Carlitz, have led to the discoveries of many interesting results (see [10,12-14,15] and references therein). Especially, the degenerate Stirling numbers of both kinds are degenerate versions of the usual Stirling numbers of both kinds and occur very frequently when we study degenerate versions of many special numbers and polynomials. In addition, $\lambda$-analogues, which are different from degenerate versions, of some special numbers and polynomials have also been studied (see \cite{5,6,9,11}). Here we consider the $\lambda$-analogues of Stirling numbers of both kinds. \par

The Stirling number of the second kind ${{n} \brace {k}}$ enumerates the number of partitions of the set $[n]=\left\{1,2,\dots,n \right\}$ into $k$ nonempty disjoint sets, while the unsigned Stirling number of the first kind ${{n} \brack {k}}=(-1)^{n-k}S_{1}(n,k)$ counts the number of permutations of $n$ elements with $k$ disjoint cycles. Here we consider the $\lambda$-analogues of the Stirling numbers of both kinds. The $\lambda$-analogues of ${{n} \brace {k}}$ are called the $\lambda$-Stirling numbers of the second kind and denoted by ${{n} \brace {k}}_{\lambda}$, while those of ${{n} \brack {k}}$ are called the $\lambda$-Stirling numbers of the first kind and denoted by ${{n} \brack {k}}_{\lambda}$. Here $\lambda$ is any real number and we note that 
${{n} \brack {k}}_{1}={{n} \brack {k}},\,\, \textrm{and} \,\,{{n} \brace {k}}_{1}={{n} \brace {k}}$. Now, let $\lambda$ be any positive integer. Then it is possible to give combinatorial interpretations of ${{n} \brace {k}}_{\lambda}$ and ${{n} \brack {k}}_{\lambda}$. Indeed, ${{n} \brace {k}}_{\lambda}$, called the translated Whitney numbers of the second kind in \cite{2,16}, enumerates the number of partitions of the set $[n]$  into $k$ subsets such that each element of each subset can mutate in $\lambda$ ways, except the dominant one. While ${{n} \brack {k}}_{\lambda}$, called the translated Whitney numbers of the first kind in \cite{2,16}, counts the number of permutations of $n$ elements with $k$ cycles such that the element of each cycle can mutate in $\lambda$ ways, except the dominant one. We also consider generalizations of the $\lambda$-Stirling numbers of both kinds, namely the $\lambda$-analogues of $r$-Stirling numbers of the second kind ${{n+r} \brace {k+r}}_{r,\lambda}$ and those of (unsigned) $r$-Stirling numbers of the first kind ${{n+r} \brack {k+r}}_{r,\lambda}$. Here we note that ${{n} \brace {k}}_{0,\lambda}={{n} \brace {k}}_{\lambda}$ and  ${{n} \brack {k}}_{0,\lambda}={{n} \brack {k}}_{\lambda}$. \par
The aim of this paper is to derive new types of generating functions of the $\lambda$-Stirling numbers of both kinds (see Theorem 3, Corollary 4, Remark 5). We note that the new type of generating function of the $\lambda$-Stirling numbers of the second kind is equivalent the ordinary generating function of those numbers (see Theorem 3, Remark 5). The generating function of the $\lambda$-Stirling numbers of the first kind follows from that of the second kind by inversion (see Corollary 4). In Theorem 6, we find expressions of the integral over $(a,\infty)$ of the reciprocals of the generalized rising factorials as finite sums by using an identity in \eqref{eq17} and the fact that the limit as $b$ tends to $\infty$ of a certain finite sum is zero (see \eqref{eq24}). Also, we derive an expression of $\frac{1}{ka^{k}}$ as an infinite sum by integrating the generating function of the $\lambda$-Stirling numbers of the second kind in Corollary 4. In Section 3, we derive the corresponding results to those ones in Section 2 for the $\lambda$-analogues of $r$-Stirling numbers of both kinds by using similar methods. We are indebted to \cite{3} for many ideas of this paper.\par

The outline of this paper is as follows. In Section 1, we recall the generalized falling factorial sequence, the Stirling numbers of the first kind, the Stirling numbers of the second kind and the unsigned Stirling numbers of the first kind. We remind the reader of the $\lambda$-Stirling numbers of the first kind, the $\lambda$-Stirling numbers of the second kind and the unsigned $\lambda$-Stirling numbers of the first kind. We recall the $\lambda$-analogues of (unsigned) $r$-Stirling numbers of the first kind and the $\lambda$-analogues of $r$-Stirling numbers of the second kind. We recall an explicit expression of the $\lambda$-analogues of $r$-Stirling numbers of the second kind and the $\lambda$-analogues of binomial coefficients. Then we remind the reader of orthogonality and inverse relations for the $\lambda$-analogues of $r$-Stirling numbers. Section2 is the main result of this paper. In Theorem 3, we obtain an expression of the reciprocal of the generalized rising factorial as an infinite series involving the $\lambda$-Stirling numbers of the second kind, which is equivalent to the ordinary generating function of those numbers. It is noted in Corollary 4 that we get, by inversion, a generating function of the unsigned $\lambda$-Stirling numbers of the first kind in terms of the reciprocals of the generalized rising factorials. By applying a Frullani's integral in \eqref{eq21}, we show that the limit as $b$ tends to $\infty$ of a certain finite sum is equal to zero (see \eqref{eq24}). By using this fact and an expression of the reciprocals of the generalized rising factorials in \eqref{eq17}, we get explicit expressions as finite sums of the integral of those reciprocals over $(a,\infty)$ in Theorem 6 and those as infinite sums of $\frac{1}{ka^{k}}$ in Corollary 7. In Section 3, all the results in Section 2 are extended to the $\lambda$-analogues of $r$-Stirling numbers and similar results to those ones in Section 2 are obtained. In the rest of this section, we recall the necessary facts that will be used throughout this paper. \\

For any $\lambda \in  \mathbb{R}$, the generalized falling factorial sequence is defined by

\begin{equation}\label{eq01}
\begin{split}
(x)_{0,\lambda}=1, \ (x)_{n,\lambda}=x(x-\lambda)(x-2\lambda) \ \cdots \ (x-(n-1)\lambda), \quad (n\geq1).
\end{split}
\end{equation}

For $n\geq0$, it is known that the Stirling numbers of the first kind are defined by
\begin{equation}\label{eq02}
\begin{split}
(x)_n=\sum_{k=0}^n S_1(n,l)x^l, \quad (\rm{see} \ \  [2-16,18]).
\end{split}
\end{equation}

As the inversion formula of \eqref{eq02}, the Stirling numbers of the second kind are defined by
\begin{equation}\label{eq03}
\begin{split}
x^n=\sum_{k=0}^n {n \brace k}(x)_k, \quad (n\geq0), \quad (\rm{see} \ \  [2-16,18,19]).
\end{split}
\end{equation}

The unsigned Stirling numbers of the first kind are defined by ${n \brack k}=(-1)^{n-k}S_1(n,k), \quad (n,k\geq0)$.
In \cite{9,11}, the $\lambda$-Stirling numbers of the first kind are given by
\begin{equation}\label{eq04}
\begin{split}
(x)_{n,\lambda}=\sum_{k=0}^n S_{1,\lambda}(n,k)x^k, \quad (n\geq0).
\end{split}
\end{equation}

The unsigned $\lambda$-Stirling numbers of the first kind are defined by
\begin{equation}\label{eq05}
\begin{split}
<x>_{n,\lambda}=(-1)^n (-x)_{n,\lambda}=\sum_{k=0}^n {n \brack k}_\lambda x^k,
\end{split}
\end{equation}
where the generalized rising factorial sequence is given by
\begin{equation} \label{eq05-1}
\begin{split}
 <x>_{0,\lambda}=1, \ <x>_{n,\lambda}=x(x+\lambda) \ \cdots \ (x+(n-1)\lambda), \quad (n\geq1).
\end{split}
\end{equation}

From \eqref{eq04} and \eqref{eq05}, we note that ${n \brack k}_{\lambda}=(-1)^{n-k}S_{1,\lambda}(n,k)$.

The $\lambda$-Stirling numbers of the second kind are defined by
\begin{equation}\label{eq06}
\begin{split}
x^n=\sum_{k=0}^n{n \brace k}_{\lambda}(x)_{k,\lambda}, \quad (n\geq0), \quad (\text{see \cite{9, 11}}).
\end{split}
\end{equation}

For $r\in\mathbb{N}\cup\{0\}$, the $\lambda$-analogues of (unsigned) $r$-Stirling numbers of the first kind are defined by
\begin{equation}\label{eq07}
\begin{split}
<x+r>_{n,\lambda}=\sum_{k=0}^n {{n+r} \brack {k+r}}_{r,\lambda} x^k, (\text{see \cite{9, 11}}),
\end{split}
\end{equation}
where we note that ${{n} \brack {k}}_{0,\lambda}={{n} \brack {k}}_{\lambda}$, for $r=0$,\quad(\text{see}\,\,\eqref{eq05}, \eqref{eq07}).

Thus, by \eqref{eq07}, we get
\begin{equation}\label{eq08}
\begin{split}
(x)_{n,\lambda}=\sum_{k=0}^n(-1)^{n-k}{{n+r} \brack {k+r}}_{r,\lambda}(x+r)^k, \quad (n\geq0).
\end{split}
\end{equation}

The $\lambda$-analogues of $r$-Stirling numbers of the second kind are defined by
\begin{equation}\label{eq09}
\begin{split}
(x+r)^n=\sum_{k=0}^n {{n+r} \brace {k+r}}_{r,\lambda}(x)_{k,\lambda}, \quad (n\geq0), \quad (\text{see \cite{9, 11}}).
\end{split}
\end{equation}

Note that
\begin{equation}\label{eq10}
\begin{split}
{{n+r} \brace {k+r}}_{r,\lambda}=\frac{1}{\lambda^k}\frac{1}{k!}\sum_{l=0}^k \binom{k}{l}(-1)^{k-l}(l\lambda+r)^n.
\end{split}
\end{equation}
This follows easily from Theorem 1 of \cite{9}. As a special case of \eqref{eq10}, we observe that 
\begin{equation}\label{eq11}
{{n} \brace {k}}_{\lambda}={{n} \brace {k}}_{0,\lambda}=\frac{1}{k!}\sum_{l=0}^k \binom{k}{l}(-1)^{k-l}\lambda^{n-k}l^{n},\quad (\text{see}\,\, \eqref{eq06}, \eqref{eq09}).
\end{equation}

For $n\geq0$, the $\lambda$-analogues of binomial coefficients are defined by
\begin{equation}\label{eq12}
\begin{split}
{\binom{x}{n}}_{\lambda}=\frac{(x)_{n,\lambda}}{n!}=\frac{x(x-\lambda) \ \cdots \ (x-(n-1)\lambda))}{n!}, \quad (\text{see \cite{15}}).
\end{split}
\end{equation}

We are going to use the next theorem and its corollary several times in the sequel. Both (a) and (b) follow from \eqref{eq08} and \eqref{eq09}, while (c), (d) and (e) can be derived from (a) and (b).

\begin{theorem}
The $\lambda$-analogues of $r$-Stirling numbers  enjoy the following orthogonality and inverse relations:
\begin{align*}
&(a)\,\, \sum_{k=l}^n (-1)^{n-k}{{n+r} \brack {k+r}}_{r,\lambda}{{k+r} \brace {l+r}}_{r,\lambda}=\delta_{n,l}, \\
&(b)\,\,\sum_{k=l}^n (-1)^{k-l}{{k+r} \brack {l+r}}_{r,\lambda}{{n+r} \brace {k+r}}_{r,\lambda}=\delta_{n,l}, \\
&(c)\,\, a_{n}=\sum_{k=0}^{n}{{n+r} \brace {k+r}}_{r,\lambda}  c_{k} \Longleftrightarrow c_{n}=\sum_{k=0}^{n}(-1)^{n-k}{{n+r}\brack {k+r}}_{r,\lambda}a_{k},\\
&(d)\,\, a_{n}=\sum_{k=n}^{m}{{k+r} \brace {n+r}}_{r,\lambda}c_{k} \Longleftrightarrow c_{n}=\sum_{k=n}^{m}(-1)^{k-n}{{k+r} \brack {n+r}}_{r,\lambda}a_{k}, \\
&(e)\,\, a_{n}=\sum_{k=n}^{\infty}{{k+r} \brace {n+r}}_{r,\lambda}c_{k} \Longleftrightarrow c_{n}=\sum_{k=n}^{\infty}(-1)^{k-n}{{k+r} \brack {n+r}}_{r,\lambda}a_{k}. \\
\end{align*}
\end{theorem}

In the special case of $r=0$, we obtain the next orthogonality and inverse relations for the $\lambda$-Stirling numbers.

\begin{corollary}
The $\lambda$-Stirling numbers  enjoy the following orthogonality and inverse relations:
\begin{align*}
&(a)\,\, \sum_{k=l}^n (-1)^{n-k}{{n} \brack {k}}_{\lambda}{{k} \brace {l}}_{\lambda}=\delta_{n,l}, \\
&(b)\,\,\sum_{k=l}^n (-1)^{k-l}{{k} \brack {l}}_{\lambda}{{n} \brace {k}}_{\lambda}=\delta_{n,l}, \\
&(c)\,\, a_{n}=\sum_{k=0}^{n}{{n} \brace {k}}_{\lambda}  c_{k} \Longleftrightarrow c_{n}=\sum_{k=0}^{n}(-1)^{n-k}{{n}\brack {k}}_{\lambda}a_{k},\\
&(d)\,\, a_{n}=\sum_{k=n}^{m}{{k} \brace {n}}_{\lambda}c_{k} \Longleftrightarrow c_{n}=\sum_{k=n}^{m}(-1)^{k-n}{{k} \brack {n}}_{\lambda}a_{k}, \\
&(e)\,\, a_{n}=\sum_{k=n}^{\infty}{{k} \brace {n}}_{\lambda}c_{k} \Longleftrightarrow c_{n}=\sum_{k=n}^{\infty}(-1)^{k-n}{{k} \brack {n}}_{\lambda}a_{k}. \\
\end{align*}
\end{corollary}

\section{New approach to $\lambda$-Stirling numbers}
Let
\begin{equation}\label{eq13}
\begin{split}
\frac{1}{x(x+\lambda)(x+2\lambda) \ \cdots \ (x+k\lambda)}=\sum_{l=0}^k \frac{A_{l,\lambda}}{x+l \lambda}.
\end{split}
\end{equation}

Then, for $0\leq m\leq k$, we have
\begin{equation}\label{eq14}
\begin{split}
\lim_{x \rightarrow -m\lambda}(x+m\lambda)\frac{1}{x(x+\lambda)(x+2\lambda) \ \cdots \ (x+k\lambda)}=\lim_{x \rightarrow -m\lambda}(x+m\lambda)\sum_{l=0}^k \frac{A_{l,\lambda}}{x+l\lambda}=A_{m,\lambda}.
\end{split}
\end{equation}

On the other hand, by \eqref{eq14}, we get
\begin{equation}\label{eq15}
\begin{split}
\lim_{x \rightarrow -m\lambda}&\frac{x+m\lambda}{x(x+\lambda) \ \cdots \ (x+(m-1)\lambda)(x+m\lambda)(x+(m+1)\lambda) \ \cdots \ (x+k\lambda)}\\
&\hspace{2cm}=\frac{1}{(-m\lambda)\cdot(-m+1)\lambda \ \cdots \ (-\lambda) \ \cdot \ \lambda \ \cdot \ 2\lambda \ \cdots \ (k-m)\lambda}\\
&\hspace{2cm}=\frac{(-1)^m}{{\lambda}^m m! \lambda^{k-m}(k-m)!}=\frac{(-1)^m}{\lambda^k}\frac{1}{k!}\binom{k}{m}.
\end{split}
\end{equation}

From \eqref{eq14} and \eqref{eq15}, we note that
\begin{equation}\label{eq16}
\begin{split}
A_{m,\lambda}=\frac{1}{\lambda^k}\frac{(-1)^m}{k!}\binom{k}{m}.
\end{split}
\end{equation}

By \eqref{eq13} and \eqref{eq16}, we get
\begin{equation}\label{eq17}
\begin{split}
\frac{1}{x(x+\lambda) \ \cdots \ (x+k\lambda)}&=\sum_{l=0}^k \frac{1}{\lambda^k}\frac{(-1)^l}{k!}\binom{k}{l}\frac{1}{x+l\lambda}\\
&=\frac{\lambda^{-k}}{k!}\sum_{l=0}^k (-1)^l \binom{k}{l}\sum_{n=0}^\infty (-1)^n l^n \lambda^n \bigg( \frac{1}{x}\bigg)^{n+1}\\
&=\sum_{n=0}^\infty \bigg(\frac{(-1)^n}{k!}\sum_{l=0}^k (-1)^l \binom{k}{l}\lambda^{n-k}l^n\bigg) \bigg(\frac{1}{x}\bigg)^{n+1}\\
&=\sum_{n=0}^\infty \frac{(-1)^{n-k}}{k!}\bigg(\sum_{l=0}^k (-1)^{k-l}\binom{k}{l}\lambda^{n-k} l^n\bigg)\bigg(\frac{1}{x}\bigg)^{n+1}\\
&=\sum_{n=k}^\infty (-1)^{n-k} {{n} \brace {k}}_\lambda \bigg(\frac{1}{x}\bigg)^{n+1},
\end{split}
\end{equation}
where we used \eqref{eq11}.

Therefore, by \eqref{eq17}, we obtain the following theorem.
\begin{theorem}
For $k\geq0$, we have
\begin{equation}\label{eq18}
\begin{split}
\frac{1}{x(x+\lambda) \ \cdots \ (x+k\lambda)}=\sum_{n=k}^\infty (-1)^{n-k} {{n} \brace {k}}_\lambda \bigg(\frac{1}{x}\bigg)^{n+1}.
\end{split}
\end{equation}
\end{theorem}

Then, from Corollary 2 (e), we obtain the following result, which is the generating function of the unsigned $\lambda$-Stirling numbers of the first kind in terms of the inverses of the generalized rising factorials.
\begin{corollary}
For $k\geq0$, we have
\begin{equation} \label{eq18-1}
\begin{split}
\bigg(\frac{1}{x}\bigg)^{k+1}&=\sum_{n=k}^\infty (-1)^{n-k}S_{1,\lambda}(n,k)\frac{1}{x(x+\lambda) \ \cdots \ (x+n\lambda)}\\
&=\sum_{n=k}^\infty {n \brack k}_{\lambda}\frac{1}{x(x+\lambda) \ \cdots \ (x+n\lambda)}.
\end{split}
\end{equation}
\end{corollary}

\medskip
\begin{remark}
It is immediate to see that \eqref{eq18} is equivalent to the ordinary generating function of the $\lambda$-Stirling numbers of the second kind given by
\begin{equation*}
\frac{x^k}{(1-\lambda x)(1-2 \lambda x) \cdots (1-k \lambda x)}=\sum_{n=k}^{\infty}{{n}\brace{k}}_{\lambda}x^{n}.
\end{equation*}
\end{remark}
On the other hand, the exponential generating function of the $\lambda$-Stirling numbers of the second kind is shown in Equation (11) of \cite{9} to be equal to
\begin{equation*}
\frac{1}{\lambda^{k}}\frac{1}{k!}\big(e^{\lambda t}-1\big)^{k}=\sum_{n=k}^{\infty}{n \brace k}_{\lambda}\frac{t^{n}}{n!}.	
\end{equation*}

Let $f$ be a nonnegative real valued function that has value at $\infty$, which is denoted by $f(\infty)$.\\
Then, for $a,b > 0$, we have
\begin{equation}\label{eq19}
\begin{split}
\int_{0}^\infty \frac{f(ax)-f(bx)}{x}dx&=\int_{0}^\infty \int_{b}^a \frac{\partial}{\partial t}f(xt)dtdx\\
&=\int_{b}^a \int_{0}^\infty \frac{\partial}{\partial t}f(xt)dxdt=\int_{b}^a \frac{f(\infty)-f(0)}{t}dt\\
&=(f(\infty)-f(0))\log\frac{a}{b}.
\end{split}
\end{equation}
The integral in \eqref{eq19} is the Frullani's integral and its evaluation in \eqref{eq20} was first published by Cauchy in 1823. The equation \eqref{eq19} is valid for `nice' functions $f$ (see \cite{1,4,17}).

Thus, by \eqref{eq19}, we get
\begin{equation}\label{eq20}
\begin{split}
\int_{0}^\infty \frac{f(ax)-f(bx)}{x}dx=(f(\infty)-f(0))\log\frac{a}{b}, \quad (\text{see \cite{4}}).
\end{split}
\end{equation}

Let us take $f(x)=e^{-x}$ and $a=1$.
Then we have
\begin{equation}\label{eq21}
\begin{split}
\int_{0}^\infty \frac{e^{-x}-e^{-bx}}{x}dx=-\log\frac{1}{b}=\log b,\quad (b >0).
\end{split}
\end{equation}

Let $k$ be any positive integer, and let $b, \lambda > 0$. Then, from \eqref{eq21}, we note that
\begin{equation}\label{eq22}
\begin{split}
\sum_{l=0}^k &\binom{k}{l}(-1)^l \log(b+l\lambda)\\
&=\sum_{l=0}^k \binom{k}{l}(-1)^l \int_{0}^\infty \frac{e^{-x}-e^{-(b+l\lambda)x}}{x}dx\\
&=\int_{0}^\infty \left\{\sum_{l=0}^k \binom{k}{l}(-1)^l \frac{e^{-x}}{x} - \sum_{l=0}^k\binom{k}{l}(-1)^l \frac{e^{-(b+l\lambda)x}}{x}\right\} dx\\
&=\int_{0}^\infty \left\{(1-1)^{k}\frac{e^{-x}}{x} - \sum_{l=0}^k\binom{k}{l}(-1)^l \frac{e^{-(b+l\lambda)x}}{x}\right\} dx\\
&=-\int_0^\infty \frac{e^{-bx}}{x}\sum_{l=0}^k \binom{k}{l}(-1)^l e^{-l\lambda x}dx=-\int_{0}^\infty \frac{e^{-bx}}{x}(1-e^{-\lambda x})^k dx.
\end{split}
\end{equation}
As $\frac{(1-e^{-\lambda x})^{k}}{x}$ is bounded on $[0,\infty)$, say by $M$, we have
\begin{equation} \label{eq23}
 \left| \int_{0}^\infty \frac{e^{-bx}}{x}(1-e^{-\lambda x})^k dx \right| \le M \int_{0}^{\infty} e^{-bx} dx =\frac{M}{b},
\end{equation}
and hence the integral $\int_{0}^\infty \frac{e^{-bx}}{x}(1-e^{-\lambda x})^k dx$ is convergent for any $b, \lambda >0$ and any positive integer $k$.

Moreover, for any positive integer $k$ and any $\lambda > 0$, by \eqref{eq23} we get
\begin{equation}\label{eq24}
\begin{split}
\lim_{b \rightarrow \infty}\sum_{l=0}^k \binom{k}{l}(-1)^l \log(b+l\lambda)=-\lim_{b \rightarrow \infty}\int_0^\infty \frac{e^{-bx}}{x}(1-e^{-\lambda x})^k dx=0.
\end{split}
\end{equation}

From \eqref{eq24} and \eqref{eq17}, for any $\lambda, a >0$ and any positive integer $k$, we note that
\begin{equation}\label{eq25}
\begin{split}
\int_a^\infty \frac{1}{x(x+\lambda) \ \cdots \ (x+k\lambda)}dx&=\lim_{b \rightarrow \infty}\int_a^b \frac{1}{x(x+\lambda) \ \cdots \ (x+k\lambda)}dx\\
&=\frac{\lambda^{-k}}{k!}\sum_{l=0}^k \binom{k}{l}(-1)^l \lim_{b \rightarrow \infty}\int_a^b \frac{1}{x+l\lambda}dx\\
&=\frac{\lambda^{-k}}{k!}\lim_{b \rightarrow \infty}\sum_{l=0}^k\binom{k}{l}(-1)^l (\log(b+l\lambda)-\log(a+l\lambda))\\
&=\frac{\lambda^{-k}}{k!}\sum_{l=0}^k\binom{k}{l}(-1)^{l+1}\log(a+l\lambda).
\end{split}
\end{equation}

Thus, by \eqref{eq24}, we get the following theorem.
\begin{theorem}
For any $\lambda, a >0$ and any positive integer $k$, the following identity holds.
\begin{equation}\label{eq26}
\begin{split}
\int_a^\infty \frac{1}{x(x+\lambda) \ \cdots \ (x+k\lambda)}dx=\frac{\lambda^{-k}}{k!}\sum_{l=0}^k\binom{k}{l}(-1)^{l+1}\log(a+l\lambda).
\end{split}
\end{equation}
\end{theorem}

From Corollary 4 and \eqref{eq26}, for any positive integer $k$ and any $\lambda, a >0$, we have
\begin{equation}\label{eq27}
\begin{split}
\frac{1}{ka^k}=\int_a^\infty \frac{1}{x^{k+1}}dx&=\sum_{n=k}^\infty {n \brack k}_\lambda \int_a^\infty \frac{1}{x(x+\lambda) \ \cdots \ (x+n\lambda)}dx\\
&=\sum_{n=k}^\infty {n \brack k}_\lambda \frac{1}{n! \lambda^n}\sum_{l=0}^n \binom{n}{l}(-1)^{l+1}\log(a+l\lambda).
\end{split}
\end{equation}

We state this as a corollary.

\begin{corollary}
For any positive integer $k$ and any $\lambda, a >0$, the following holds true.
\begin{equation*}
\frac{1}{ka^k}=\sum_{n=k}^\infty {n \brack k}_\lambda \frac{1}{n! \lambda^n}\sum_{l=0}^n \binom{n}{l}(-1)^{l+1}\log(a+l\lambda).
\end{equation*}
\end{corollary}

By Theorem 3, we get
\begin{equation}\label{eq28}
\begin{split}
\bigg(\frac{1}{x}\bigg)^k= \sum_{n=k}^\infty {n \brack k}_\lambda \frac{1}{(x+\lambda)(x+2\lambda)\ \cdots \ (x+n\lambda)}.
\end{split}
\end{equation}

Replacing $x$ by $x-\lambda$ and $n$ by $n+1$, we get
\begin{equation}\label{eq29}
\begin{split}
\bigg(\frac{1}{x-\lambda}\bigg)^k =\sum_{n=k-1}^\infty {{n+1} \brack k}_\lambda \frac{1}{x(x+\lambda) \ \cdots \ (x+n\lambda)}.
\end{split}
\end{equation}

In particular, for $k=1$, we have
\begin{equation}\label{eq30}
\begin{split}
\frac{1}{x-\lambda}&=\sum_{n=0}^{\infty}\frac{\lambda^{n}n!}{x(x+\lambda) \ \cdots \ (x+n\lambda)}=\sum_{n=0}^\infty \frac{\lambda^{n}}{\binom{x+(n-1)\lambda}{n}_\lambda (x+n\lambda)}.
\end{split}
\end{equation}

\bigskip

\section{$\lambda$-analogues of $r$-Stirling numbers}

Let
\begin{equation}\label{eq31}
\begin{split}
\frac{1}{(x+r)(x+r+\lambda) \ \cdots \ (x+r+k\lambda)}=\sum_{l=0}^k \frac{K_{l,\lambda}}{x+r+l\lambda}.
\end{split}
\end{equation}

For $0\leq m\leq k$, we have
\begin{equation}\label{eq32}
\begin{split}
K_{m,\lambda}&=\lim_{x \rightarrow {-r-m\lambda}}\sum_{l=0}^k \frac{K_{l,\lambda}}{x+r+l\lambda}(x+r+m\lambda)\\
&=\lim_{x \rightarrow {-r-m\lambda}}\frac{x+r+m\lambda}{(x+r)(x+r+\lambda) \ \cdots \ (x+r+k\lambda)}\\
&=\frac{(-1)^m}{\lambda^k}\frac{k!}{m!(k-m)!}\frac{1}{k!}=\frac{(-1)^m}{\lambda^k}\frac{1}{k!}\binom{k}{m}.
\end{split}
\end{equation}

By \eqref{eq31} and \eqref{eq32}, we get
\begin{equation}\label{eq33}
\begin{split}
\frac{1}{(x+r)(x+r+\lambda) \ \cdots \ (x+r+k\lambda)}&=\sum_{l=0}^k \frac{(-1)^l}{\lambda^k k!}\binom{k}{l}\frac{1}{x+r+l\lambda}\\
&=\frac{1}{\lambda^k k!}\sum_{l=0}^k(-1)^l\binom{k}{l}\frac{1}{x}  \frac{1}{1+\frac{r+l\lambda}{x}}\\
&=\sum_{n=0}^\infty (-1)^n \frac{1}{\lambda^k k!}\sum_{l=0}^k(-1)^l\binom{k}{l}(r+l\lambda)^n \bigg(\frac{1}{x}\bigg)^{n+1}\\
&=\sum_{n=0}^\infty (-1)^{n-k}\bigg(\frac{1}{\lambda^k}\frac{1}{k!}\sum_{l=0}^k(-1)^{k-l}\binom{k}{l}(r+l\lambda)^n\bigg)\bigg(\frac{1}{x}\bigg)^{n+1}\\
&=\sum_{n=k}^\infty (-1)^{n-k}{{n+r} \brace {k+r}}_{r,\lambda}\bigg(\frac{1}{x}\bigg)^{n+1},
\end{split}
\end{equation}
where we used \eqref{eq10}.

Therefore, by \eqref{eq33}, we obtain the following theorem.

\begin{theorem}
For $k\geq0$, we have
\begin{equation}\label{eq34}
\begin{split}
\frac{1}{(x+r)(x+r+\lambda) \ \cdots \ (x+r+k\lambda)}=\sum_{n=k}^\infty (-1)^{n-k}{{n+r} \brace {k+r}}_{r,\lambda} \bigg(\frac{1}{x}\bigg)^{n+1}.
\end{split}
\end{equation}
\end{theorem}

\medskip

Now, from Corollary 2 (e) and using the notation in \eqref{eq12}, we obtain the following result..
\begin{corollary}
For $k\geq0$, we have
\begin{equation*}
\begin{split}
\bigg(\frac{1}{x}\bigg)^{k+1}&=\sum_{n=k}^\infty {{n+r} \brack {k+r}}_{r,\lambda} \frac{1}{(x+r)(x+r+\lambda) \ \cdots \ (x+r+n\lambda)}\\
&=\sum_{n=k}^\infty {{n+r} \brack {k+r}}_{r,\lambda}\frac{1}{n!{\binom{x+r+(n-1)\lambda}{n}}_\lambda (x+r+n\lambda)}.
\end{split}
\end{equation*}
\end{corollary}

\medskip

\begin{remark}
The equation \eqref{eq34} is equivalent to the ordinary generating function of the $\lambda$-analogues of $r$-Stirling numbers of the second kind given by
\begin{equation*}
\frac{x^{k}}{(1-rx)(1-(r+\lambda)x)\cdots(1-(r+k \lambda)x)}=\sum_{n=k}^{\infty}{{n+r} \brace {k+r}}_{r,\lambda}x^{n}.
\end{equation*}
On the other hand, it is shown in Theorem 1 of \cite{9} that the exponential generating function of the $\lambda$-analogues of $r$-Stirling numbers of the second kind is given by
\begin{equation*}
\frac{1}{\lambda^{k}}\frac{1}{k!}\big(e^{\lambda t}-1\big)^{k}e^{rt}=\sum_{n=k}^{\infty}{n+r \brace k+r}_{r,\lambda}\frac{t^{n}}{n!},\quad (k\ge 0). 
\end{equation*}
\end{remark}
Let $k$ be any positive integer and let $\lambda, a >0$.
By using \eqref{eq24} and \eqref{eq33}, we obtain
\begin{equation}\label{eq35}
\begin{split}
\int_a^\infty \frac{1}{(x+r)(x+r+\lambda) \ \cdots \ (x+r+k\lambda)}dx=\frac{1}{k!}\frac{1}{\lambda^k}\sum_{l=0}^k \binom{k}{l}(-1)^{l+1}\log(a+r+l\lambda).
\end{split}
\end{equation}

From \eqref{eq35} and Corollary 6, we have
\begin{equation}\label{eq36}
\begin{split}
\frac{1}{ka^k}&=\int_a^\infty \frac{1}{x^{k+1}}dx\\
&=\sum_{n=k}^\infty {{n+r} \brack {k+r}}_{r,\lambda}\frac{1}{n!} \int_a^\infty \frac{1}{n!{\binom{x+r+(n-1)\lambda}{n}}_\lambda (x+r+n\lambda)}dx\\
&=\sum_{n=k}^\infty {{n+r} \brack {k+r}}_{r,\lambda}\frac{1}{n!\lambda^n}\sum_{l=0}^n \binom{n}{l}(-1)^{l+1}\log(a+r+l\lambda).
\end{split}
\end{equation}

Thus we have shown the following theorem.
\begin{theorem}
For any positive integer $k$ and any $\lambda, a >0$, the following hold true.
\begin{align*}
&\int_a^\infty \frac{1}{(x+r)(x+r+\lambda) \ \cdots \ (x+r+k\lambda)}dx=\frac{1}{k!}\frac{1}{\lambda^k}\sum_{l=0}^k \binom{k}{l}(-1)^{l+1}\log(a+r+l\lambda), \\
&\frac{1}{ka^k}=\sum_{n=k}^\infty {{n+r} \brack {k+r}}_{r,\lambda}\frac{1}{n!\lambda^n}\sum_{l=0}^n \binom{n}{l}(-1)^{l+1}\log(a+r+l\lambda).
\end{align*}
\end{theorem}

\vspace{1cm}

\section{Conclusion}
We have witnessed that studying both degenerate versions and $\lambda$-analogues of some special numbers and polynomials yielded many fascinating results. Here we note that these two are different in nature. Indeed, the degenerate versions tend to the non-degenerate ones as $\lambda$ tends to 0, while the $\lambda$-analogues as $\lambda$ tends to 1. In this paper, we studied the $\lambda$-Stirling numbers of both kinds which are $\lambda$-analogues of the usual Stirling numbers. We derived new types of generating functions for those numbers and also for the $\lambda$-analogues of $r$-Stirling numbers of both kinds. All of these generating functions are related to reciprocals of the generalized rising factorials. Furthermore, from those new types of generating functions we obtained expressions of the integral over $(a,\infty)$ of the reciprocals of the generalized rising factorials as finite sums and an expression of $\frac{1}{ka^{k}}$ as an infinite sum involving the $\lambda$-Stirling numbers of the first kind or the $\lambda$-analogues of $r$-Stirling numbers of the first kind. \par
It is one of our future projects to continue to study various $\lambda$-analogues of many special numbers and polynomials and to find their applications to physics, science and engineering as well as to mathematics.

\end{document}